\documentclass[reqno, a4paper]{amsart}

\usepackage{amsmath,amsfonts,amssymb}
\usepackage{verbatim}
\usepackage{enumerate}
\usepackage{url}
\usepackage{tikz}
\usepackage[left=2cm,right=2cm, top=3cm, bottom=2cm]{geometry}
\usepackage{slashbox}

\usetikzlibrary{arrows}

\def\N{\mathbb N}

\def\Z{\mathbb Z}

\def\ord{\mathop{\rm ord}\nolimits}

\theoremstyle{plain}
\newtheorem{theorem}{Theorem}[section]
\newtheorem{lemma}[theorem]{Lemma}

\def\qed{\hfill\hbox{$\square$}}

\theoremstyle{definition}

\newtheorem*{Ack}{Acknowledgements}
\numberwithin{equation}{section}

\title{Extremal product-one free sequences in Dihedral and Dicyclic groups}

\author[F. E. Brochero Mart\'{\i}nez]{F. E. Brochero Mart\'{\i}nez}
\address{
Departamento de Matem\'{a}tica\\
Universidade Federal de Minas Gerais\\
UFMG\\
Belo Horizonte, MG\\
30123-970\\
Brazil\\
}
\email{fbrocher@mat.ufmg.br }

\author[S\'avio Ribas]{S\'avio Ribas}
\email{savio.ribas@gmail.com }

\date{\today}

\begin{document}

\maketitle

\begin{abstract}
Let $G$ be a finite group, written multiplicatively. The Davenport constant of $G$ is the smallest positive integer $D(G)$ such that every sequence of $G$ with $D(G)$ elements has a non-empty subsequence with product $1$. Let $D_{2n}$ be the Dihedral Group of order $2n$ and $Q_{4n}$ be the Dicyclic Group of order $4n$. J. J. Zhuang and W. Gao \cite{ZhGa} showed that $D(D_{2n}) = n+1$ and J. Bass \cite{Bas} showed that $D(Q_{4n}) = 2n+1$. In this paper, we give explicit characterizations of all sequences $S$ of $G$ such that $|S| = D(G) - 1$ and $S$ is free of subsequences whose product is 1, where $G$ is equal to $D_{2n}$ or $Q_{4n}$ for some $n$. 
\end{abstract}

\section{Introduction}

Given a finite group $G$ written multiplicatively, the {\em Zero-Sum Problems} study conditions that guarantee a non-empty subsequence in $G$ has a non-empty subsequence with certain prescribed properties such that the product of its elements, in some order, is equal to the identity $1 \in G$. Interesting properties include length, repetitions, and weights. 

One of the first problem of this type was considered by Erd\"os, Ginzburg, and Ziv \cite{EGZ}. They proved that given $2n-1$ integers, it is possible to select $n$ of them, such that their sum is divisible by $n$. In the language of group theory, every sequence $S$ consisting of at least $2n-1$ elements in a finite cyclic group of order $n$ has a subsequence $T$ of length $n$ such that the product of the elements in $T$ in some order is equal to the identity. Further, they proved that the number $2n-1$ is the smallest positive integer with this property. 

These problems have been studied extensively for abelian groups; see the surveys by Y. Caro \cite{Car} and W. Gao - A. Geroldinger \cite{GaGe}.

An important problem of zero-sum type is to determine the so-called {\em Davenport constant} of a finite group $G$ (written multiplicatively). This constant, denoted by $D(G)$, is the smallest positive integer $d$ such that every sequence $S$ with $d$ elements in $G$ (repetition allowed) contains some subsequence $T$ such that the product of the elements in $T$ in some order is $1$.

For $n \in \N$, let $C_n \simeq \Z_n$ denote the cyclic group of order $n$ written multiplicatively. The Davenport constant is known for the following groups:
\begin{itemize}
\item $D(C_n) = n$;
\item $D(C_m \times C_n) = m + n - 1$ if $m|n$  (J. Olson, \cite{Ols2});
\item $D(C_{p^{e_1}} \times \dots \times C_{p^{e_r}}) = 1 + \sum_{i = 1}^r (p^{e_i} - 1)$ (J. Olson, \cite{Ols1});
\item $D(D_{2n}) = n + 1$ where $D_{2n}$ is the Dihedral Group of order $2n$ (see \cite{ZhGa});
\item $D(C_q \rtimes_s C_m) = m + q - 1$ where $q \ge 3$ is a prime number and $\ord_q(s) = m \ge 2$ (J. Bass, \cite{Bas}).
\end{itemize}
For most finite groups, the Davenport constant is not known. % (what about upper and lower bounds?).

By the definition of the Davenport constant, for a given finite group $G$ there exist sequences $S$ with elements in $G$ such that $|S| = D(G) - 1$ and such that $S$ is \emph{free of product-1 subsequences}. That is, there exists $S = (g_1, \cdots, g_{D(G) - 1})$ of $G$ such that $g_{i_1} \cdots g_{i_k} \ne 1$ for every non-empty subset $\{i_1, \cdots, i_k\}$ of $\{1, \cdots, D(G) - 1\}$. The {\em Inverse Zero-Sum Problems} study the structure of these extremal-length sequences which are free of product-$1$ subsequences with some prescribed property. For an overview, see articles by 
W. Gao, A. Geroldinger, D. J. Grynkiewicz and W. A. Schmid \cite{GaGeSc}, \cite{Sch} and \cite{GaGeGr}.

The inverse problems associated with Davenport constant are solved for few abelian groups. For example, the following theorem is direct consequence of Theorem \ref{propertyC} and gives a complete characterization of sequences free of product-$1$ subsequences in the finite cyclic group $C_n$. 

\begin{theorem}\label{propC}
Let $S$ be a sequence in $C_n$ free of product-$1$ subsequences with $n-1$ elements, where $n \ge 2$. Then
\begin{equation*}
S = (\underbrace{g,\dots,g}_{n-1 \text{ times}}), \text{ where $g$ is a generator of $C_n$.}
\end{equation*}
\end{theorem}

%Observe that in the cyclic case, sequences free of product-1 subsequences are exactly the diagonal sequences. In particular, such sequences have a single element repeated many times.
Observe that in the cyclic case, sequences free of product-$1$ subsequences contain an element repeated many times. It is natural to ask if this is true in general. Specifically, let us say that a given finite abelian group $G$ has {\em Property C} if every maximal sequence $S$ of $G$ free of product-$1$ subsequences with at most $\exp(G)$, the exponent of the group $G$, elements has the form 
$$S = (\underbrace{T, T, \dots, T}_{\exp(G) - 1 \text{ times}}),$$
for some subsequence $T$ of $S$. The above theorem states that $C_n$ has Property C. It then follows from a result of C. Reiher \cite{Rei} that $C_p^2$ possesses Property C (see also \cite{GaGeSc} and \cite{GaGe3}). In \cite{GaGeGr}, W. Gao, A. Geroldinger and D. J. Grynkiewicz showed that this result is multiplicative, extending our conclusion to the groups $C_n^2$ for a composite number $n$. In \cite{Sch2}, W. A. Schmid discusses the case $C_n \times C_m$, where $n|m$. Not much is known about groups of rank exceeding $2$, only few specific cases (see, for example, \cite{Sch2}).

A {\em minimal zero sequence} $S$ in a finite abelian group $G$ is a sequence such that the product of its elements is $1$, but each proper subsequence is free of product-$1$ subsequences. In \cite[Theorem~6.4]{GaGe2}, W. Gao and A. Geroldinger showed that if $|S| = D(G)$ then $S$ contains some element $g \in G$ with order $\ord(g) = \exp(G)$ for certain groups such as $p$-groups, cyclic groups, groups with rank two and groups that are the sum of two elementary $p$-groups. They also conjectured that the same conclusion holds for every finite abelian group.

For non-abelian groups, nothing was known, until the paper \cite{BR} where we solved the inverse problem associated to the metacyclic group
$$C_q \rtimes_s C_m = \langle x,y| x^m=1, y^q=1, yx = xy^s, \ord_q(s)=m, q \text{ prime} \rangle.$$

Specifically, we proved the following result:

\begin{theorem}[{\cite[Theorem~1.2]{BR}}] \label{metacyclic}
Let $q$ be a prime number, $m \ge 2$ be a divisor of $q-1$ and $s \in \Z_q^*$ such that $\ord_q(s) = m$. 
%Let $q$ be a prime number, $m \geq 2$ a divisor of $q-1$ and $s$ a non-zero element of $(\Z/q\Z)$ such that $\ord_q(s) = m$. (The reason is that it is now usually understood that $\Z_p$ refers to the $p$-adic integers) 

Let $S$ be a sequence in the metacyclic group $C_q \rtimes_s C_m$ with $m+q-2$ elements.
\begin{enumerate}
\item If $(m,q) \neq (2,3)$ then the following statements are equivalent:
\begin{enumerate}[(i)]
\item $S$ is free of product-$1$ subsequences;
\item For some $1 \le t \le q-1$, $1 \le i \le m-1$ such that $\gcd(i,m) = 1$ and $0 \le \nu_1, \dots, \nu_{m-1} \le q-1$,
$$S = ( \underbrace{y^t, y^t, \dots, y^t}_{q-1\text{ times}}, x^iy^{\nu_1}, x^iy^{\nu_2}, \dots, x^iy^{\nu_{m-1}} ).$$
\end{enumerate}
\item If $(m,q) = (2,3)$ then $S$ is free of product-$1$ subsequences if and only if 
\begin{equation*}
S = (y^t,y^t,xy^{\nu}) \text{ for }t \in \{2,3\} \text{ and } \nu \in \{0,1,2\} \;\;\; \text{ or } \;\;\; S = (x,xy,xy^2).
\end{equation*}
\end{enumerate}
\end{theorem}

In this article, we characterize the maximal sequences which are free of product-$1$ subsequences for the {\em Dihedral Groups} and the {\em Dicyclic Groups}. In particular, we show that these sequences have a property similar to Property C.

Let $n \ge 2$. Denote by $D_{2n} \simeq C_n \rtimes_{-1} C_2$ the {\em Dihedral Group} of order $2n$, i.e., the group generated by $x$ and $y$ with relations:
\begin{equation*}\label{defdihedral}
x^2 = y^n = 1, \quad yx = xy^{-1}.
\end{equation*}

Let $n \ge 2$. Denote by $Q_{4n}$ the {\em Dicyclic Group}, i.e., the group generated by $x$ and $y$ with relations:
\begin{equation*}\label{defdicyclic}
x^2 = y^n, \quad y^{2n} = 1, \quad yx = xy^{-1}.
\end{equation*}

We have $Z(Q_{4n}) = \{1,y^n\}$, where $Z(G)$ denotes the center of a group $G$. In addition: 
\begin{equation*}
Q_{4n}/\{1,y^n\} \simeq D_{2n},
\end{equation*}
where $D_{2n}$ is the Dihedral Group of order $2n$.

\vspace{0.2cm}

Specifically, we prove the following result:

\begin{theorem}\label{inversedihedral}
Let $S$ be a sequence in the Dihedral Group $D_{2n}$ with $n$ elements, where $n \ge 3$.
\begin{enumerate}
\item If $n \ge 4$ then the following statements are equivalent:
\begin{enumerate}[(i)]
\item $S$ is free of product-$1$ subsequences;
\item For some $1 \le t \le n-1$ with $\gcd(t,n)=1$ and $0 \le s \le n-1$,
$$S = ( \underbrace{y^t, y^t, \dots, y^t}_{n-1\text{ times}}, xy^s ).$$
\end{enumerate}
\item If $n=3$ then $S$ is free of product-$1$ subsequences if and only if 
\begin{equation*}
S = (y^t,y^t,xy^{\nu}) \text{ for }t \in \{2,3\} \text{ and } \nu \in \{0,1,2\} \;\;\; \text{ or } \;\;\; S = (x,xy,xy^2).
\end{equation*}
\end{enumerate}
\end{theorem}

Notice that this theorem reduces to Theorem \ref{metacyclic} in the case when $n$ is prime and when $m=2$, since we have $s\equiv-1 \pmod n$. For $n \ge 4$, it is easy to check that $(ii) \Longrightarrow (i)$. Observe that the case $n = 3$ is exactly the item $(2)$ of Theorem \ref{metacyclic}, since $D_6 \simeq C_3 \rtimes_{-1} C_2$. For $n=2$, we have $D_4 \simeq \Z_2^2$ and it is easy to check that $S$ is free of product-$1$ subsequences if and only if $S = (x,y)$, $S = (xy,y)$ or $S = (x,xy)$.

\vspace{0.2cm}

As a consequence of previous theorem, we obtain:

\begin{theorem}\label{inversedicyclic}
Let $S$ be a sequence in the Dicyclic Group $Q_{4n}$ with $2n$ elements, where $n \ge 2$.
\begin{enumerate}
\item If $n \ge 3$ and $|S| = 2n$ then the following statements are equivalent:
\begin{enumerate}[(i)]
\item $S$ is free of product-$1$ subsequences;
\item For some $1 \le t \le n-1$ with $\gcd(t,2n) = 1$ and $0 \le s \le 2n-1$,
$$S = ( \underbrace{y^t, y^t, \dots, y^t}_{2n-1\text{ times}}, xy^s ).$$
\end{enumerate}
\item If $n = 2$ then $S$ is free of product-$1$ subsequences if and only if, for some $r \in \Z_4^*$ and $s \in \Z_4$, $S$ has one of the forms
$$(y^r,y^r,y^r,xy^s), (y^r,xy^s,xy^s,xy^s) \text{ or }(xy^s,xy^s,xy^s,xy^{r+s}).$$
\end{enumerate}
\end{theorem}

Again, if $n \ge 3$ then it is easy to check that $(ii) \Longrightarrow (i)$, therefore we just need to show that $(i) \Longrightarrow (ii)$. If $n = 2$ then it is easy to check that the sequences of the form $(y^r,y^r,y^r,xy^s)$, $(y^r,xy^s,xy^s,xy^s)$, $(xy^s,xy^s,xy^s,xy^{r+s})$ are free of product-$1$ subsequences, therefore we just need to show that all other sequences with $2n$ elements have subsequences with product $1$.

Also note that if $n = 2$ then $Q_8$ is isomorphic to the {\em Quaternion Group}, i.e. the group defined by
$$\langle e,i,j,k| \; i^2 = j^2 = k^2 = ijk = e, e^2 = 1 \rangle.$$
This isomorphism may be described, for example, by $x \mapsto i, y \mapsto j$ (or its natural permutations, by the simmetry of $i,j,k$). The above theorem says that, in terms of Quarternion Group, extremal sequences free of product-$1$ subsequence are of the forms
$$\pm(i,i,i,\pm j), \pm(i,i,i,\pm k), \pm(j,j,j,\pm i), \pm(j,j,j,\pm k), \pm(k,k,k,\pm i) \text{ or } \pm(k,k,k,\pm j).$$

The main technical difficulty in our present proofs lies on the fact that $n$ may not be prime, therefore we cannot use Cauchy-Davenport inequality (see \cite[p. 44-45]{Nat}) and Vosper's Theorem (see \cite{Vos}), as used in the proof of Theorem \ref{metacyclic} (see \cite{BR}). Instead, we now exhibit sequences with product-$1$ in cases not covered by those given forms. In section \ref{proofdihedral} we prove Theorem \ref{inversedihedral}, solving the extremal inverse zero-sum problem associated to Davenport constant for dihedral groups. The proof of Theorem \ref{inversedicyclic} is split up into sections \ref{proofdicyclic}, \ref{proofdicyclic3} and \ref{proofdicyclic2}, where we solve the extremal inverse zero-sum problem associated to Davenport constant for dicyclic groups of orders $4n$, where $n \ge 4$, $n = 3$ and $n = 2$, respectively. The case $n \ge 4$ is a direct consequence of Theorem \ref{inversedihedral} and the case $n = 3$ follows from Theorem \ref{inversedihedral}, but in the special case $(2)$. The proof in the case $n = 2$ is done manually, using only Theorem \ref{propertyC} to reduce the number of cases, without using Theorem \ref{inversedihedral}.

\vspace{0.2cm}

\section{Notation and Auxiliary Results}

In this section we present the notation and auxiliary lemmas and theorems that we use throughout the paper.

Let $G$ be a finite group written multiplicatively and $S = (g_1, g_2, \dots, g_l)$ be a sequence of elements of $G$. We denote by $|S| = l$ the length of $S$.

For each subsequence $T$ of $S$, i.e. $T = (g_{n_1}, g_{n_2}, \dots, g_{n_k})$, where $\{n_1, n_2, \dots, n_k\}$ is a subset of $\{1, 2, \dots, l\}$, we say that $T$ is a {\em product-$1$ subsequence} when
$$g_{\sigma(n_1)}g_{\sigma(n_2)}\dots g_{\sigma(n_k)} = 1$$
 for some permutation $\sigma$ of $\{n_1,\dots,n_k\}$, and if there are no such product-$1$ subsequences then we say that $S$ is {\em free of product-$1$ subsequences}.

If $S_1 = (g_{i_1},\dots, g_{i_u})$ and $S_2 = (g_{j_1}, \dots, g_{j_v})$ are subsequences of $S$, then
\begin{itemize}
\item $SS_1^{-1}$ denotes the subsequence formed by the elements of $S$ without the elements of $S_1$;
\item $S_1S_2 = (g_{i_1}, \dots, g_{i_u}, g_{j_1}, \dots, g_{j_v})$ denotes the concatenation of $S_1$ and $S_2$;
\item $S^k = SS \dots S$ denotes the concatenation of $k$ identical copies of $S$'s.
\end{itemize}

For the group $G = D_{2n} = \langle x,y| x^2=1,\  y^n=1,\  yx=xy^{-1} \rangle$, let
\begin{itemize}
\item $H_D$ be the normal cyclic subgroup of order $n$ generated by $y$;
\item $N_D = D_{2n} \setminus H_D = x \cdot H_D$.
\end{itemize}
The product of any even number of elements in $N_D$ is in $H_D$, since
\begin{equation}\label{produtoDN}
xy^{\alpha} \cdot xy^{\beta} = y^{\beta - \alpha}.
\end{equation}

For the group $G = Q_{4n} = \langle x,y| x^2=y^n,\  y^{2n}=1,\  yx=xy^{-1} \rangle$, let
\begin{itemize}
\item $H_Q$ be the normal cyclic subgroup of order $2n$ generated by $y$;
\item $N_Q = Q_{4n} \setminus H_Q = x \cdot H_Q$.
\end{itemize}
The product of any even number of elements in $N_Q$ is in $H_Q$, since
\begin{equation}\label{produtoQN}
xy^{\alpha} \cdot xy^{\beta} = y^{\beta - \alpha + n}.
\end{equation}

\vspace{0.2cm}

The following theorem is known as ``Davenport constant of $\mathbb Z_n$ with weights $\{\pm 1\}$'', and will be used in the proof of Theorem \ref{inversedihedral}.

\begin{lemma}[{\cite[Lemma~2.1]{ACFKP}}] \label{weighted}
Let $n \in \mathbb N$ and $(y_1,\dots,y_s)$ be a sequence of integers with $s > \log_2 n$. Then there exist a nonempty $J \subset \{1,2,3,\dots,s\}$ and $\varepsilon_j \in \{\pm 1\}$ for each $j \in J$ such that
\begin{equation*}
\sum_{j \in J} \varepsilon_j y_j \equiv 0 \pmod n.
\end{equation*}
\end{lemma}

As a generalization of Theorem \ref{propC}, we will use the following:

\begin{theorem}[{\cite[Theorem~2.1]{GaGeSc}, see also \cite{BEN}}] \label{propertyC}
Let $G$ be a cyclic group of order $n \ge 3$ and $S$ be a zero-sum free sequence in $G$ of length $|S| \ge (n+1)/2$. Then there exists some $g \in S$ that appear at least $2|S| - n + 1$ times in $S$. In particular, $D(G) = n$ and the following statements hold:
\begin{enumerate}
\item If $|S| = n-1$, then $S = (g)^{n-1}$.
\item If $|S| = n-2$, then either $S = (g)^{n-2}$ or $S = (g)^{n-3}(g^2)$.
\item If $|S| = n-3$, then $S$ has one of the following forms: $$(g)^{n-3}, (g)^{n-4}(g^2), (g)^{n-4}(g^3), (g)^{n-5}(g^2)^2.$$
\end{enumerate}
\end{theorem}

\vspace{0.2cm}

\section{Proof of Theorem \ref{inversedihedral}}\label{proofdihedral}

We just need to show that $(i) \Longrightarrow (ii)$. Let $S$ be a sequence in $D_{2n}$ with $n$ elements that is free of product-$1$ subsequences. If $S \cap N_D$ contains two identical elements then $S$ is not free of product-$1$ subsequences by Equation \ref{produtoDN}. Hence, we assume that the elements of $S \cap N_D$ are all distinct.

From now on, we consider some cases, depending on the cardinality of $S \cap H_D$:

\vspace{0.2cm}

\begin{enumerate}[(a)]
\item {\bf Case $|S \cap H_D| = n$:} In this case, $S$ is contained in the cyclic subgroup of order $n$. Since $D(\mathbb Z_n) = n$, $S$ contains some non-empty subsequence with product $1$.

\vspace{0.3cm}

\item {\bf Case $|S \cap H_D| = n-1$:} In this case, by Theorem \ref{propertyC}, the elements of $S \cap H_D$ must all be equal, say, $S \cap H_D = (y^t)^{n-1}$ where $\gcd(t,n) = 1$, and so $S = (y^t)^{n-1}(xy^s)$.

\vspace{0.3cm}

\item {\bf Case $|S \cap H_D| = n-2$:} In this case, by Theorem \ref{propertyC}, $S$ must have one of these forms:

\vspace{0.2cm}

\begin{enumerate}[(c-1)]
\item {\bf Subcase $S = (y^t)^{n-2}(xy^u,xy^v)$:} Notice that we can obtain the products $$xy^u \cdot \underbrace{y^t \cdot y^t \dots y^t}_{k \text{ times}} \cdot xy^v = y^{v-u-kt}.$$
Since $1 \le k \le n-2$, it's enough to take $k \equiv (v-u)t^{-1} \pmod n$. The only problem occurs when $k = n-1$, that is, when $v - u + t \equiv 0 \pmod n$, but in this case we switch $u$ and $v$ and so $$xy^v \cdot y^t \cdot xy^u = y^{u-v-t} = 1.$$

\vspace{0.2cm}

\item {\bf Subcase $S = (y^t)^{n-3}(y^{2t},xy^u,xy^v)$:} Notice that we can obtain the products $$xy^u \cdot \underbrace{y^t \cdot y^t \dots y^t}_{k \text{ times}} \cdot xy^v = y^{v-u-kt}.$$
Since $1 \le k \le n-1$, it's enough to take $k \equiv (v-u)t^{-1} \pmod n$.
\end{enumerate}

\vspace{0.3cm}

\item {\bf Case $|S \cap H_D| = n-3$:} 
In this case, by Theorem \ref{propertyC}, $S$ must contain at least $n-5$ copies of some $y^t$, where $\gcd(t,n)=1$. Also, suppose that
$$S \cap N_D = (xy^{\alpha}, xy^{\beta}, xy^{\gamma}).$$
By renaming $z = y^t$, we may assume without loss of generality that $t = 1$. By Pigeonhole Principle it follows that there exist two exponents of $y$ with difference in $\{1, 2, \dots, \lfloor n/3 \rfloor\} \pmod n$, say,
$$\alpha - \beta \in \{0, 1, 2, \dots, \lfloor n/3 \rfloor\} \pmod n.$$
We may ensure that
\begin{equation*}
xy^{\beta} \cdot y^r \cdot xy^{\alpha} = y^{\alpha - \beta - r} = 1
\end{equation*}
for some $1 \le r \le n/3$ provided there are enough $y$'s. But if $n \ge 8$ then $r \le n/3 \le n-5$, so there are enough $y$'s and the theorem follows in these cases. Since the theorem is already proved for $n$ prime, it only remains to prove for $n \in \{4,6\}$.

\vspace{0.3cm}

For $n = 4$, we have $S = (xy^{\alpha},xy^{\beta},xy^{\gamma},y)$ and some of the $\alpha,\beta,\gamma$ are consecutives modulo $4$. Without loss of generality, suppose that $\alpha - \beta \equiv 1 \pmod 4$, so $xy^{\beta}\cdot y \cdot xy^{\alpha} = 1$. 

\vspace{0.3cm}

For $n = 6$, we have the cases
\begin{equation*}
\begin{cases}
S = (xy^{\alpha},xy^{\beta},xy^{\gamma},y,y,y); \\
S = (xy^{\alpha},xy^{\beta},xy^{\gamma},y,y,y^2); \\
S = (xy^{\alpha},xy^{\beta},xy^{\gamma},y,y,y^3); \text{ or } \\
S = (xy^{\alpha},xy^{\beta},xy^{\gamma},y,y^2,y^2).
\end{cases}
\end{equation*}
If the set $\{\alpha,\beta,\gamma\} \pmod 6$ contains two consecutive elements, say $\alpha \equiv \beta + 1 \pmod 6$, then
$$xy^{\alpha} \cdot xy^{\beta} \cdot y = y^{\beta + 1 - \alpha} = 1.$$
Otherwise, the only possibilities are $\{\alpha,\beta,\gamma\} \pmod 6 = \{0,2,4\}$ or $\{\alpha,\beta,\gamma\} \pmod 6 = \{1,3,5\}$. Suppose that $\alpha \equiv \beta + 2 \pmod 6$ and notice that in any option for $S$ it is possible to take a product $y^2$ coming from $S \cap H_D$. Therefore
$$xy^{\alpha} \cdot xy^{\beta} \cdot y^2 = y^{\beta + 2 - \alpha} = 1.$$

\vspace{0.3cm}

\item {\bf Case $|S \cap H_D| = n-k, 4 \le k \le n$:} Suppose that 
\begin{align*}
S \cap H_D &= (y^{t_1}, y^{t_2}, \dots, y^{t_{n-k}}), \\
S \cap N_D &= (xy^{\alpha_1}, xy^{\alpha_2}, \dots, xy^{\alpha_k}).
\end{align*}
It follows from Theorem \ref{weighted} that if
\begin{equation}\label{desigualdadehipotese}
\lfloor k/2 \rfloor > \left\lfloor \log_2 n \right\rfloor
\end{equation}
then there exist a linear combination of a subset of
\begin{equation*}
\{ (\alpha_1 - \alpha_2), \; (\alpha_3 - \alpha_4), \; \dots, \; (\alpha_{2\lfloor k/2 \rfloor - 1} - \alpha_{2\lfloor k/2 \rfloor}) \}
\end{equation*}
with coefficients $\pm 1$ summing $0$. Suppose without loss of generality that, in this combination,
\begin{equation*}
(\alpha_1 - \alpha_2), \; (\alpha_3 - \alpha_4), \; \dots, \; (\alpha_{2u - 1} - \alpha_{2u})
\end{equation*}
appear with signal $-1$ and
\begin{equation*}
(\alpha_{2u + 1} - \alpha_{2u + 2}), \; (\alpha_{2u + 3} - \alpha_{2u + 4}), \; \dots, \; (\alpha_{2v - 1} - \alpha_{2v}),
\end{equation*}
appear with signal $+1$. Then
\begin{align*}
(xy^{\alpha_1} \cdot xy^{\alpha_2}) \dots (xy^{\alpha_{2u-1}} \cdot xy^{\alpha_{2u}}) \cdot (xy^{\alpha_{2u+2}} \cdot xy^{\alpha_{2u+1}}) \cdot (xy^{\alpha_{2v}} \cdot xy^{\alpha_{2v - 1}}) &= \\
y^{\alpha_2 - \alpha_1} \dots y^{\alpha_{2u} - \alpha_{2u-1}} \cdot y^{\alpha_{2u+1} - \alpha_{2u+2}} \dots y^{\alpha_{2v - 1} - \alpha_{2v}} &= 1
\end{align*}
Thus the theorem is true for $k > 2 \lfloor \log_2 n \rfloor$ + 1.

Hence, we may assume $4 \le k \le 2\lfloor \log_2 n \rfloor + 1$. Theorem \ref{propertyC} implies without loss of generality that $t_i = t$ for $1 \le i \le n-2k+1$. By renaming $z = y^t$, we may assume without loss of generality that $t = 1$. Since $k \ge 4$, Pigeonhole Principle implies that there exist $\alpha_i, \alpha_j$ such that 
\begin{equation*}
\alpha_i - \alpha_j \in \{1, 2, \dots, \lfloor n/4 \rfloor\}.
\end{equation*}
Notice that if 
\begin{equation}\label{equationnk}
n - 2k + 1 \ge n/4
\end{equation}
then
$$xy^{\alpha_j} \cdot y^r \cdot xy^{\alpha_i} = y^{\alpha_i - \alpha_j - r} = 1$$
for some $0 \le r \le \lfloor n/4 \rfloor$. But if $n \ge 8$ then $3n \ge 8 \log_2 n \ge 4(k - 1)$, therefore Equation \ref{equationnk} holds in this case. Since the theorem is already proved for $n$ prime, it only remains to prove for $n \in \{4,6\}$.

\vspace{0.3cm}

For $n = 4$, the only possibility is $S = (x,xy,xy^2,xy^3)$. Thus, $x \cdot xy \cdot xy^3 \cdot xy^2 = 1$.

\vspace{0.3cm}

For $n = 6$, there are three subcases to consider:

\vspace{0.2cm}

\begin{itemize}
\item {\bf Subcase $k = 4$:} Let $S = (y^{t_1}, y^{t_2}, xy^{\alpha_1}, xy^{\alpha_2}, xy^{\alpha_3}, xy^{\alpha_4})$. Then either there are two pairs of consecutive $\alpha_i$'s modulo $6$ or there are three consecutive $\alpha_i$'s modulo $6$. 

\begin{itemize}
\item
If there are two pairs of consecutive $\alpha_i$'s, say $\alpha_1 + 1 \equiv \alpha_2$ and $\alpha_3 + 1 \equiv \alpha_4 \pmod 6$, then
$$xy^{\alpha_1} \cdot xy^{\alpha_2} \cdot xy^{\alpha_4} \cdot xy^{\alpha_3} = 1.$$

\item
If there are three consecutive $\alpha_i$'s, say $\alpha_1 + 2 \equiv \alpha_2 + 1 \equiv \alpha_3 \pmod 6$, then $\alpha_4$ can be any element in the set $\{\alpha_3 + 1, \alpha_3 + 2, \alpha_3 + 3\}$. If $\alpha_4 \equiv \alpha_3 + 1$ or $\alpha_4 \equiv \alpha_3 + 3 \pmod 6$ then we return to the previous item. Therefore we may assume $\alpha_4 \equiv \alpha_3 + 2 \pmod 6$. Taking the products $xy^{\alpha_i} \cdot xy^{\alpha_j} = y^{\alpha_j - \alpha_i}$, we can get any element in $\{y, y^2, y^3, y^4, y^5\}$. For example,
\begin{align*}
xy^{\alpha_1} \cdot xy^{\alpha_2} &= y, \\
xy^{\alpha_1} \cdot xy^{\alpha_3} &= y^2, \\
xy^{\alpha_2} \cdot xy^{\alpha_4} &= y^3, \\
xy^{\alpha_1} \cdot xy^{\alpha_4} &= y^4, \text{ and } \\
xy^{\alpha_2} \cdot xy^{\alpha_1} &= y^5.
\end{align*}
Let $i,j$ such that $\alpha_i - \alpha_j \equiv t_1 \pmod 6$. Then $$xy^{\alpha_i} \cdot xy^{\alpha_j} \cdot y^{t_1} = y^{\alpha_j - \alpha_i + t_1} = 1.$$
\end{itemize}

\vspace{0.2cm}

\item {\bf Subcase $k = 5$:} Let $S = (y^t, xy^{\alpha_1}, xy^{\alpha_2}, \dots, xy^{\alpha_5})$. In this case, there are four consecutive $\alpha_i$'s modulo $6$, say $\alpha_1 + 3 \equiv \alpha_2 + 2 \equiv \alpha_3 + 1 \equiv \alpha_4 \pmod 6$. Then
$$xy^{\alpha_1} \cdot xy^{\alpha_2} \cdot xy^{\alpha_3} \cdot xy^{\alpha_4} = y^{\alpha_2 - \alpha_1 + \alpha_4 - \alpha_3} = 1.$$

\vspace{0.2cm}

\item {\bf Subcase $k = 6$:} The only possibility is $S = (x,xy,xy^2,xy^3,xy^4,xy^5)$. So
$$x \cdot xy \cdot xy^3 \cdot xy^2 = 1.$$
\end{itemize}
\end{enumerate}
\qed

\section{Proof of Theorem \ref{inversedicyclic} in the case $n \ge 4$}\label{proofdicyclic}

We just need to prove that $(i) \Longrightarrow (ii)$. Let $S$ be a sequence in $Q_{4n}$ with $2n$ elements that is free of product-$1$ subsequences. We consider some cases depending on the cardinality of $S \cap H_Q$:

\vspace{0.2cm}

\begin{enumerate}[({5.}1)]
\item {\bf Case $|S \cap H_Q| = 2n$:} In this case, $S$ is contained in the cyclic subgroup of order $2n$. Since $D(H_Q) = D(\Z_{2n}) = 2n$, $S$ contains some non-empty subsequence with product $1$.

\vspace{0.3cm}

\item {\bf Case $|S \cap H_Q| = 2n-1$:} In this case, by Theorem \ref{propertyC}, the elements of $S \cap H_Q$ must all be equal, say, $S \cap H_Q = (y^t)^{2n-1}$ where $\gcd(t,2n) = 1$, and so $S = (y^t)^{2n-1}(xy^s)$. 

\vspace{0.3cm}

\item {\bf Case $|S \cap H_Q| = 2n-2$:} In this case, let $S_1$ be a subsequence of $S$ such that $|S_1 \cap H_Q| = n-2$ and let $S_2 = S S_1^{-1}$. Then $S_2$ is a sequence in $H_Q$ with $n$ elements. Since $$Q_{4n} / \{1,y^n\} \simeq D_{2n},$$ Theorem \ref{inversedihedral} tells us that $S_1$ and $S_2$ must contain subsequences $T_1$ and $T_2$, respectively, with products in $\{1,y^n\}$. If some of these products is $1$ then we are done. Otherwise, both products are $y^n$, therefore $$\prod_{z \in T_1T_2} z = y^{2n} = 1,$$ thus $S$ is not free of product-$1$ subsequences.

\vspace{0.3cm}

\item {\bf Case $|S \cap H_Q| = 2n-3$:} In this case, let $S_1$ be a subsequence of $S$ such that $|S_1 \cap H_Q| = n-3$ and let $S_2 = S S_1^{-1}$. Then $S_2$ is a sequence in $H_Q$ with $n$ elements. The argument is similar to the above case, thus $S$ is not free of product-$1$ subsequence.

\vspace{0.3cm}

\item {\bf Case $|S \cap H_Q| = 2n-k, 4 \le k \le 2n$:} In this case, let $S_1$ be a subsequence of $S$ such that $|S_1 \cap H_Q| \le n-2$ and $S_2 = S S_1^{-1}$ is such that $|S_2 \cap H_Q| \le n-2$. The argument is similar to the above cases, thus $S$ is not free of product-$1$ subsequence.
\end{enumerate}

\vspace{0.3cm}

Therefore, the proof for $n \ge 4$ is complete.
\qed

\vspace{0.3cm}

\section{Proof of Theorem \ref{inversedicyclic} in the case $n = 3$}\label{proofdicyclic3}

We have $Q_{12} = \langle x,y| \; x^2 = y^3, y^6 = 1, yx = xy^5 \rangle$ and we just need to prove that $(i) \Longrightarrow (ii)$. Let $S$ be a sequence in $Q_{12}$ with $6$ elements that is free of product-$1$ subsequences. We consider some cases depending on the cardinality of $S \cap H_Q$:

\vspace{0.2cm}

\begin{enumerate}[({6.}1)]
\item {\bf Case $|S \cap H_Q| = 6$:} In this case, $S$ must contain a product-$1$ subsequence, since $D(H_Q) = D(\Z_6) = 6$.

\vspace{0.3cm}

\item {\bf Case $|S \cap H_Q| = 5$:} In this case, Theorem \ref{propertyC} says that $S \cap H_Q = (y^r)^5$ where $r \in \{1,5\}$, therefore $S$ is of the form $(y^r)^5(xy^s)$.

\vspace{0.3cm}

\item {\bf Case $|S \cap H_Q| = 4$:} In this case, we decompose $S = S_1S_2$ where $|S_i| = 3$ for $i \in \{1,2\}$, $|S_1 \cap H_Q| = 1$ and $|S_2 \cap H_Q| = 3$, and use the same argument than item $(5.3)$, therefore $S$ is not free of product-$1$ subsequences.

\vspace{0.3cm}

\item {\bf Case $|S \cap H_Q| = 3$:} In this case, we decompose $S = S_1S_2$ where $|S_i| = 3$ for $i \in \{1,2\}$, $|S_1 \cap H_Q| = 1$ and $|S_2 \cap H_Q| = 2$, therefore
\begin{align*}
S_1 \!\!\! \pmod {\{1,y^3\}} &= (y^r,xy^u,xy^v) \text{ and } \\
S_2 \!\!\! \pmod {\{1,y^3\}} &= (y^t,y^t,xy^s) \text{ for } r,t \in \{1,2\} \text{ and } s,u,v \in \{0,1,2\}.
\end{align*}
Notice that $S_1$ contains a subsequence with product in $\{1,y^3\}$. Observe that if $r \neq t$ then we could also decompose $S = S_1'S_2'$ where 
$$S_1' = S_1(y^t)(y^r)^{-1} \;\;\;\;\; \text{ and } \;\;\;\;\; S_2' = S_2(y^r)(y^t)^{-1}.$$
So, $S_1'$ and $S_2'$ have subsequences with products in $\{1,y^3\}$ and we can use the same argument than item $(5.3)$. Therefore, $r = t$ and we can decompose $S = S_1''S_2''$ such that
$$S_1'' \!\!\! \pmod {\{1,y^3\}} = (y^t)^3 \;\;\;\;\; \text{ and } \;\;\;\;\; S_2'' \!\!\! \pmod {\{1,y^3\}} = (xy^s,xy^u,xy^v).$$
Notice that $S_1''$ contains a subsequence with product in $\{1,y^3\}$ and $S_2''$ does not contains a subsequence with product in $\{1,y^3\}$ if and only if $S_2'' \! \pmod {\{1,y^3\}} = (x,xy,xy^2)$. Hence, the only possibility for $S \! \pmod {\{1,y^3\}}$ is:
$$S \!\!\! \pmod {\{1,y^3\}} = (y^t)^3(x,xy,xy^2),$$
and so the possibilities for $S \cap H_Q$ are
$$ (y)^3, \;\; (y^2)^3, \;\; (y^4)^3, \;\; (y^5)^3, \;\; (y)^2(y^4), \;\; (y)(y^4)^2, \;\; (y^2)^2(y^5) \;\; \text{ or } \;\; (y^2)(y^5)^2.$$
Notice that the second, third, fifth and seventh possibilities contain subsequences with product $1$, therefore it only remains
$$(y)^3, \;\; (y^5)^3, \;\; (y)(y^4)^2 \;\; \text{ or } \;\; (y^2)(y^5)^2.$$
Observe that in every case, it is possible to find a subsequence with product either $y$ or $y^5$.
We claim that $S \cap N_Q$ contains subsequences with product $y$ and $y^5$, which we can join with those $y$ or $y^5$ coming from $S \cap H_Q$ to get a product-$1$ subsequence. For this, a sufficient condition is the existence of two elements in $S \cap N_Q$ such that the exponents of $y$ have difference $2$, since $xy^{\alpha} \cdot xy^{\alpha+2} = y$ and $xy^{\alpha+2} \cdot xy^{\alpha} = y^5$. In fact, if $xy^{\beta} \in S \cap N_Q$ and $xy^{\beta-2},xy^{\beta+2} \not\in S \cap N_Q$ then $xy^{\beta+1},xy^{\beta-1} \in S \cap N_Q$, and so $(\beta+1) - (\beta-1) = 2$.

\vspace{0.3cm}

\item {\bf Case $|S \cap H_Q| = 2$:} In this case, we decompose $S = S_1S_2$ where $|S_i| = 3$ and $|S_i \cap H_Q| = 1$, and use the same argument than item $(5.3)$, therefore $S$ is not free of product-$1$ subsequences.

\vspace{0.3cm}

\item {\bf Case $|S \cap H_Q| = 1$:} In this case, we decompose $S = S_1S_2$ where $|S_i| = 3$, $|S_1 \cap H_Q| = 1$ and $|S_2 \cap H_Q| = 0$. Notice that $S_2$ contains a subsequence with product in $\{1,y^3\}$. The same argument than item $(5.3)$ does not apply if and only if $S_2 \pmod {\{1,y^3\}} = (x,xy,xy^2)$. Observe that we could also decompose $S = S_1'S_2'$ in such way that $|S_i'| = 3$ for $i \in \{1,2\}$, $|S_1' \cap H_Q| = 1$, $|S_2' \cap H_Q| = 0$ and $S_2' \pmod {\{1,y^3\}}$ contains two elements with the same exponent in $y$. Therefore, the same argument than item $(5.3)$ applies for $S_1'$ and $S_2'$.

\vspace{0.3cm}

\item {\bf Case $|S \cap H_Q| = 0$:} In this case, we decompose $S = S_1S_2$ where $|S_i| = 3$ and $|S_i \cap H_Q| = 0$ for $i \in \{1,2\}$. Notice that $S_i$ contains a subsequence with product in $\{1,y^3\}$ if and only if $S_i \pmod {\{1,y^3\}}$ is not of the form $(x,xy,xy^2)$. Hence, the same argument than item $(5.3)$ does not apply if and only if $S_i \pmod {\{1,y^3\}} = (x,xy,xy^2)$ for $i \in \{1,2\}$. Observe that we could also decompose $S = S_1'S_2'$ in such way that $S_1' \pmod {\{1,y^3\}} = (x,x,xy)$ and $S_2' \pmod {\{1,y^3\}} = (xy,xy^2,xy^2)$. Therefore, the same argument than item $(5.3)$ applies for $S_1'$ and $S_2'$.
\end{enumerate}

\vspace{0.3cm}

Therefore, the proof for $n = 3$ is complete.
\qed

\vspace{0.3cm}

\section{Proof of Theorem \ref{inversedicyclic} in the case $n = 2$}\label{proofdicyclic2}

We have $Q_8 = \langle x,y| \; x^2 = y^2, y^4 = 1, yx = xy^3 \rangle$. Suppose that $S$ is a sequence in $Q_8$ with $4$ elements that is free of product-$1$ subsequences. We want to show that $S$ has some of those forms given in item $(2)$. For this, we consider some cases depending on the cardinality of $S \cap H_Q$:

\vspace{0.3cm}

\begin{enumerate}[({7.}1)]
\item {\bf Case $|S \cap H_Q| = 4$:} In this case, $S$ must contain a product-$1$ subsequence, since $D(H_Q) = D(\Z_4) = 4$.

\vspace{0.3cm}

\item {\bf Case $|S \cap H_Q| = 3$:} In this case, Theorem \ref{propertyC} says that $S \cap H_Q = (y^r)^3$ where $r \in \{1,3\}$, therefore $S$ is of the form $(y^r)^3(xy^s)$.

\vspace{0.3cm}

\item {\bf Case $|S \cap H_Q| = 2$:} In this case, the only possibilities for $S \cap H_Q$ making $S$ be free of product-$1$ subsequences are $$(y)^2, (y,y^2), (y^3)^2 \;\;\; \text{ and } \;\;\; (y^2,y^3),$$ and in all these possibilities $S \cap H_Q$ possesses a subsequence with product $y^2$, namely $y \cdot y$, $y^3 \cdot y^3$ or $y^2$ itself. On the other hand, the possibilities for $S \cap N_Q$ are $$(x,xy), (x,xy^2), (x,xy^3), (xy,xy^2), (xy,xy^3), (xy^2,xy^3) \;\;\; \text{ and } \;\;\; (xy^s,xy^s)$$ for $s \in \Z_4$. \\
The cases $(xy^s,xy^s)$ can be eliminated, since $xy^s \cdot xy^s \cdot y^2 = 1$. \\
The cases $(x,xy^2)$ and $(xy,xy^3)$ can also be eliminated, since $xy^s \cdot xy^{s+2} = 1$. \\
The other cases can be eliminated by the following table:

\vspace{0.3cm}

\begin{center}
\begin{tabular}{c||c|c|c|c}
\backslashbox{$S \cap N_Q$}{$S \cap H_Q$}	&	$(y)^2$	&	$(y,y^2)$	&	$(y^3)^2$	&	$(y^2,y^3)$ \\ \hline \hline
$(x,xy)$		&	$x \cdot xy \cdot y = 1$	&	$x \cdot xy \cdot y = 1$	&	$xy \cdot x \cdot y^3 = 1$	&	$xy \cdot x \cdot y^3 = 1$ \\ \hline
$(x,xy^3)$		&	$x \cdot y \cdot xy^3 = 1$	&	$x \cdot y \cdot xy^3 = 1$	&	$x \cdot xy^3 \cdot y^3 = 1$	&	$x \cdot xy^3 \cdot y^3 = 1$ \\ \hline
$(xy,xy^2)$		&	$xy \cdot xy^2 \cdot y = 1$	&	$xy \cdot xy^2 \cdot y = 1$	&	$xy^2 \cdot xy \cdot y^3 = 1$	&	$xy^2 \cdot xy \cdot y^3 = 1$ \\ \hline
$(xy^2,xy^3)$	&	$xy^2 \cdot xy^3 \cdot y = 1$	&	$xy^2 \cdot xy^3 \cdot y = 1$	&	$xy^3 \cdot xy^2 \cdot y^3 = 1$	&	$xy^3 \cdot xy^2 \cdot y^3 = 1$ 
\end{tabular}
\end{center}

\vspace{0.3cm}
 
\item {\bf Case $|S \cap H_Q| = 1$:} In this case, the possibilities for $S \cap H_Q$ are $(y)$, $(y^2)$ and $(y^3)$. On the other hand, $S \cap N_Q$ has three elements and, by Equation \ref{produtoQN}, we may assume that $x \in S$. Thus, the possibilities for $S \cap N_Q$ are

\begin{center}
$(x,x,x)$, $(x,x,xy)$, $(x,x,xy^2)$, $(x,x,xy^3)$, $(x,xy,xy)$, $(x,xy,xy^2)$, $(x,xy,xy^3),$

\vspace{0.1cm}

$(x,xy^2,xy^2)$, $(x,xy^2,xy^3)$ and $(x,xy^3,xy^3).$
\end{center}

\vspace{0.3cm}

\noindent If $S$ contains $(x,xy^2)$ or $(xy,xy^3)$ then $S$ is not free of product-$1$ subsequences, since
\begin{equation}\label{Qdiferenca2}
x \cdot xy^2 = 1 = xy \cdot xy^3.
\end{equation}
Therefore, the remainder possibilities are $(x,x,x)$, $(x,x,xy)$, $(x,x,xy^3)$, $(x,xy,xy)$ and $(x,xy^3,xy^3)$. Notice that these last four possibilities contains two identical terms and two terms such that the exponents of $y$ have difference $1$ modulo $4$. Since 
\begin{align*}
xy^{\alpha} \cdot xy^{\alpha+1} \cdot y &= 1 \\
xy^{\alpha} \cdot xy^{\alpha} \cdot y^2 &= 1 \\
xy^{\alpha+1} \cdot xy^{\alpha} \cdot y^3 &= 1,
\end{align*}
we may discard these four cases. Therefore, the only remainder possibility is $(x,x,x)$, and so $S = (y)(xy^s)^3$ or $S = (y^3)(xy^s)^3$.

\vspace{0.3cm}

\item {\bf Case $|S \cap H_Q| = 0$:} In this case, we also may assume $x \in S$. Therefore, the possibilities for $S$ are

\begin{center}
$(x,x,x,x)$, $(x,x,x,xy)$, $(x,x,x,xy^2)$, $(x,x,x,xy^3)$, $(x,x,xy,xy)$, $(x,x,xy,xy^2),$

\vspace{0.1cm}

$(x,x,xy,xy^3)$, $(x,x,xy^2,xy^2)$, $(x,x,xy^2,xy^3)$, $(x,x,xy^3,xy^3)$, $(x,xy,xy,xy),$

\vspace{0.1cm}

$(x,xy,xy,xy^2)$, $(x,xy,xy,xy^3)$, $(x,xy,xy^2,xy^2)$, $(x,xy,xy^2,xy^3),$

\vspace{0.1cm}

$(x,xy^2,xy^2,xy^2)$, $(x,xy^2,xy^2,xy^3)$, $(x,xy^2,xy^3,xy^3)$ and $(x,xy^3,xy^3,xy^3).$
\end{center}
\vspace{0.2cm}

\noindent If $S$ contains two pairs of identical elements, say $(x,x,x^{\alpha},x^{\alpha})$, then
$$x \cdot x \cdot xy^{\alpha} \cdot xy^{\alpha} = 1,$$
so we remove $(x,x,x,x)$, $(x,x,xy,xy)$, $(x,x,xy^2,xy^2)$ and $(x,x,xy^3,xy^3)$. \\
If $S$ contains some of the pairs $(x,xy^2)$ or $(xy,xy^3)$ then, by equations in \ref{Qdiferenca2}, we may remove other $11$ possibilities, thus it only remains $(x,x,x,xy)$, $(x,x,x,xy^3)$, $(x,xy,xy,xy)$ and $(x,xy^3,xy^3,xy^3)$. Therefore, $S = (xy^s)(xy^{r+s})^3$ for $r \in \Z_4^*$ and $s \in \Z_4$.
\end{enumerate}
\qed

\Ack The second author would like to thank {\em CAPES}/Brazil for the PhD student fellowship.

\end{document}